\documentclass[a4paper,12pt]{amsart}
\usepackage{amssymb}
\usepackage{ifthen}
 \usepackage[dvips]{graphicx}
\nonstopmode \numberwithin{equation}{section}
\usepackage{amssymb}
\usepackage{ifthen}
\usepackage{graphicx}
\usepackage{amsmath}
\usepackage[T1]{fontenc} 
\usepackage[utf8]{inputenc}
\usepackage[usenames,dvipsnames]{color}
\usepackage{color}
\usepackage[english]{babel}
\usepackage{fancyhdr}
\usepackage{fancybox}
\usepackage{tikz}

\nonstopmode \numberwithin{equation}{section}
\theoremstyle{plain}
\newtheorem{thm}[equation]{Theorem}
\newtheorem{cor}[equation]{Corollary}
\newtheorem{lem}[equation]{Lemma}
\newtheorem{prop}{Proposition}

\newtheorem{conj}{Conjecture}

\theoremstyle{definition}
\newtheorem{defn}{Definition}[section]

\newtheorem{prob}{Problem}
\newtheorem{rem}{Remark}[section]

\newcounter{minutes}
\divide\time by 60
\newcounter{hours}
\multiply\time by 60
\addtocounter{minutes}{-\time}

\newcounter {own}
\def\theown {\thesection       .\arabic{own}}

\newenvironment{pf}[1][]{%
 \vskip 3mm
 \noindent
 \ifthenelse{\equal{#1}{}}%
  {{\slshape Proof. }}%
  {{\slshape #1.} }%
 }%
{\qed\bigskip}

\newcounter{alphabet}

\newcommand{\real}{{\operatorname{Re}\,}}

\def\be{\begin{equation}}
\def\ee{\end{equation}}

\newcommand{\bee}{\begin{enumerate}}
\newcommand{\eee}{\end{enumerate}}

\newcommand{\blem}{\begin{lem}}
\newcommand{\elem}{\end{lem}}
\newcommand{\bthm}{\begin{thm}}
\newcommand{\ethm}{\end{thm}}
\newcommand{\bcor}{\begin{cor}}
\newcommand{\ecor}{\end{cor}}
\newcommand{\beg}{\begin{examp}}
\newcommand{\eeg}{\end{examp}}
\newcommand{\begs}{\begin{examples}}
\newcommand{\eegs}{\end{examples}}

\newcommand{\bdefn}{\begin{defn}}
\newcommand{\edefn}{\end{defn}}

\newcommand{\bprob}{\begin{prob}}
\newcommand{\eprob}{\end{prob}}
\newcommand{\bei}{\begin{itemize}}
\newcommand{\eei}{\end{itemize}}

\newcommand{\bcon}{\begin{conj}}
\newcommand{\econ}{\end{conj}}
\newcommand{\bcons}{\begin{conjs}}
\newcommand{\econs}{\end{conjs}}
\newcommand{\bprop}{\begin{prop}}
\newcommand{\eprop}{\end{prop}}
\newcommand{\br}{\begin{rem}}
\newcommand{\er}{\end{rem}}
\newcommand{\brs}{\begin{rems}}
\newcommand{\ers}{\end{rems}}
\newcommand{\bo}{\begin{obser}}
\newcommand{\eo}{\end{obser}}
\newcommand{\bos}{\begin{obsers}}
\newcommand{\eos}{\end{obsers}}
\newcommand{\bpf}{\begin{pf}}
\newcommand{\epf}{\end{pf}}
\newcommand{\ba}{\begin{array}}
\newcommand{\ea}{\end{array}}
\newcommand{\beq}{\begin{eqnarray}}
\newcommand{\beqq}{\begin{eqnarray*}}
\newcommand{\eeq}{\end{eqnarray}}
\newcommand{\eeqq}{\end{eqnarray*}}

\begin{document}

\title{Bohr-Rogosinski and improved Bohr type inequalities for certain fully starlike harmonic mappings}

\author{Molla Basir Ahamed}
\address{Molla Basir Ahamed,
	School of Basic Sciences,
	Indian Institute of Technology Bhubaneswar,
	Bhubaneswar-752050, Odisha, India.}
\email{mba15@iitbbs.ac.in}

\author{Vasudevarao Allu}
\address{Vasudevarao Allu,
School of Basic Sciences,
Indian Institute of Technology Bhubaneswar,
Bhubaneswar-752050, Odisha, India.}
\email{avrao@iitbbs.ac.in}

\subjclass[{AMS} Subject Classification:]{Primary 30C45, 30C50, 30C80}

\def\thefootnote{}
\footnotetext{ {\tiny File:~\jobname.tex,
printed: \number\year-\number\month-\number\day,
          \thehours.\ifnum\theminutes<10{0}\fi\theminutes }
} \makeatletter\def\thefootnote{\@arabic\c@footnote}\makeatother

\begin{abstract}
The classical Bohr inequality states that if $ f $ is an analytic function with the power series representation $ f(z)=\sum_{n=0}^{\infty}a_nz^n $ in the unit disk $ \mathbb{D}:=\{z\in\mathbb{C} : |z|<1\} $ such that $ |f(z)|\leq 1 $ for all $ z\in\mathbb{D} $, then 
\begin{equation*}
	\sum_{n=0}^{\infty}|a_n|r^n\leq 1\;\; \text{for}\;\; |z|=r\leq\frac{1}{3} 
\end{equation*}
and the constant $ 1/3 $ cannot be improved. The constant $ r_0=1/3 $ is known as Bohr radius and the inequality $ \sum_{n=0}^{\infty}|a_n|r^n\leq 1 $ is known as Bohr inequality. Let $ \mathcal{H} $ be the class of complex-valued harmonic mappings $ f=h+\bar{g}$ defined in the unit disk $ \mathbb{D} $, where $ h $ and $ g $ are analytic functions in $ \mathbb{D} $ with the normalization $ h(0)=0=h^{\prime}(0)-1 $ and $ g(0)=0 $. 
Let $ \mathcal{H}_{0}=\{f=h+\bar{g}\in\mathcal{H} : g^{\prime}(0)=0\}. $  Let  
$ \mathcal{P}^{0}_{\mathcal{H}}(M) :=\{f=h+\overline{g} \in \mathcal{H}_{0}: \real (zh^{\prime\prime}(z))> -M+|zg^{\prime\prime}(z)|,\; z \in \mathbb{D},\; M>0\} $. Functions in the class $ \mathcal{P}^{0}_{\mathcal{H}}(M) $ are called fully starlike univalent functions for $ 0<M<1/\log 4 $. In this paper, we obtain the sharp Bohr-Rogosinski type inequality and improved Bohr inequality and the corresponding Bohr radius for the class $ \mathcal{P}_{\mathcal{H}}^{0}(M) $.
\end{abstract}

\maketitle
\pagestyle{myheadings}
\markboth{Molla Basir Ahamed and Vasudevarao Allu}{Bohr-Rogosinski and improved Bohr type inequalities for certain fully starlike harmonic mappings}
\section{Introduction}
 Let $ f $ be an analytic function with the power series representation $ f(z)=\sum_{n=0}^{\infty}a_nz^n $ in the unit disk $ \mathbb{D}:=\{z\in\mathbb{C} : |z|<1\} $ such that $ |f(z)|\leq 1 $ for all $ z\in\mathbb{D} $. The Bohr inequality is given by 
\begin{equation}\label{e-1.1}
	M_f(r):=\sum_{n=0}^{\infty}|a_n|r^n\leq 1,\;\; \text{for all}\;\; |z|=r\leq\frac{1}{3} 
\end{equation}
and the constant $ 1/3 $ cannot be improved. The inequality \eqref{e-1.1} was established by Bohr \cite{Bohr-1914} in $ 1914 $. The constant $ r_0=1/3 $ is known as Bohr radius, while the inequality $ \sum_{n=0}^{\infty}|a_n|r^n\leq 1 $ is known as Bohr inequality. Bohr actually obtained the inequality \eqref{e-1.1} for $ r\leq 1/6 $ and subsequently later, Weiner, Riesz and Schur have independently obtained the inequality \eqref{e-1.1} for the radius $ 1/3 $, which is sharp.  Moreover, for 
\begin{equation*}
	\phi_b(z)=\frac{b-z}{1-bz},\;\; b\in [0,1),
\end{equation*}
it follows easily that $ M_{\phi_b}(r)>1 $ if, and only if, $ r>1/(1+2b) $. Consequently, the radius $ 1/3 $ is optimal when $ b\rightarrow 1 $.\vspace{1.5mm} \par For $ f(z)=\sum_{n=0}^{\infty}a_nz^n $ in $ \mathbb{D} $, the Bohr inequality can be written in terms of the Euclidean distance as follows
\begin{equation}\label{e-11.2}
	d\left(\sum_{n=0}^{\infty}|a_nz^n|,|a_0|\right)=\sum_{n=1}^{\infty}|a_nz^n|\leq 1-|f(0)|=d(f(0),\partial\mathbb{D}),
\end{equation}
where $ d $ is the Euclidean distance and $ \partial\mathbb{D} $ is the boundary of the unit disk $ \mathbb{D}. $ More generally, a class $ \mathcal{F} $ of analytic (or harmonic) functions $ f(z)=\sum_{n=0}^{\infty}a_nz^n $ mapping $ \mathbb{D} $ into a domain $ \Omega $ is said to satisfy a Bohr phenomenon if an inequality of type \eqref{e-11.2} hols uniformly in $ |z|\leq \rho_0 $, where $ 0<\rho_0\leq 1 $ for functions in the class $ \mathcal{F} $.  In $ 2004 $, Beneteau \emph{et al.} \cite{Bene-2004} first introduced the notion of the Bohr phenonomenon for a Banach space of analytic functions in the unit disk $ \mathbb{D} $. Aizenberg \emph{et al.} \cite{Aizenberg-PAMS-2000} have presented an abstract apporach to the the problem of finding Bohr radius and proved that Bohr's phenomenon occurs in a general condition.\vspace{1.5mm}

\par In $ 1995 $, Dixon \cite{Dixon & BLMS & 1995} showed via methods related to Bohr's inequality for power series \cite{Bohr-1914}, that there are Banach algebras $ \mathcal{A} $ which are not isometrically operator algebras and yet satisfy a special case of von Neumann inequality: $ ||p(a)||\leq \sup\{|p(z)| : z\in\mathbb{D}\} $ for all $ a\in\mathcal{A} $ and polynomial $ p(z) $ with $ p(0)=0 $. In the recent years, the generalization of Bohr's inequality for various classes of analytic functions becomes an active research area (see \cite{kayumov-2018-b,Kayumov-CRACAD-2020,Liu-2020}). The Bohr radius, for the holomorphic functions has been studied by Aizenberg \textit{et al.} \cite{aizenberg-2001},  Aytuna and Djakov \cite{Ayt & Dja & BLMS & 2013}. The Bohr phenomenon for the class of starlike log-harmonic mappings has been studied by Ali \textit{et al.} \cite{Ali & Abdul & Ng & CVEE & 2016}, and Ahamed and Allu \cite{Ahamed-Allu-2021}. In $ 2018 $, Kayumov and Ponnusamy \cite{kayumov-2018-b} introduced the notion of $ p $-Bohr radius for harmonic functions, and established result obtaining $ p $-Bohr radius for the class of odd analytic functions. The improved Bohr inequality for locally univalent harmonic mappings has been sudied by Evdoridis \emph{et al.} \cite{Evdoridis-IM-2018}.  In $ 2021 $, Huang \emph{et al.} \cite{Huang-Liu-Ponnusamy-MJM-2021} determined Bohr inequality for the class of harmonic mappings $ f=h+\overline{g} $ in the unit disk $ \mathbb{D} $, where both $ h(z)=\sum_{n=0}^{\infty}a_{pn+m}z^{pn+m} $ and $ g(z)=\sum_{n=0}^{\infty}b_{pn+m}z^{pn+m} $ are analytic and bounded in $ \mathbb{D} $ and also investigated Bohr-type inequalities of harmonic mappings with a multiple zero at the origin. For more exciting aspects of Bohr phenomenon, we refer to \cite{Abu-2010,Abu-2014,Aizn-PAMS-2000,Ali-2017,Alkhaleefah-PAMS-2019,Evdoridis-IM-2018,Nirupam-CVEE-2018,Ghosh-Vasudevarao-BAMS-2020,Kay & Pon & AASFM & 2019,Kayumov-CRACAD-2020,Kay & Pon & Sha & MN & 2018,Liu-Ponnusamy-BMMS-2019} and the references therein. \vspace{2mm}

\par  In addition to the Bohr inequality, the concept of the Rogosinski inequality is also used in \cite{Rogosinski-1923}, which is defined as follows: Let $f(z)=\sum_{n=0}^{\infty} a_{n}z^{n}$ be analytic in $\mathbb{D}$ and its corresponding partial sum of $f$ is defined by $S_{N}(z):=\sum_{n=0}^{N-1} a_{n}z^{n}$. Then, for every $N \geq 1$, we have $|\sum_{n=0}^{N-1} a_{n}z^{n}|<1$ in the disk $|z|<1/2$ and the radius $1/2$ is sharp. Motivated by Rogosinski radius for bounded analytic functions in $\mathbb{D}$, Kayumov and Ponnusamy \cite{kayumov-2017} have introduced Bohr-Rogosinski radius and considered the Bohr-Rogosinski sum $R_{N}^{f}(z)$ which is defined by 
\begin{equation*}
	R_{N}^{f}(z):=|f(z)|+ \sum_{n=N}^{\infty} |a_{n}||z|^{n}.
\end{equation*}
\noindent It is easy to see that $|S_{N}(z)|=\big|f(z)-\sum_{n=N}^{\infty} a_{n}z^{n}\big| \leq |R_{N}^{f}(z)|$. Thus, the validity of Bohr-type radius for $R_{N}^{f}(z)$, which is related to the classical Bohr sum (Majorant series) in which $f(0)$ is replaced by $f(z)$, gives Rogosinski radius in the case of bounded analytic functions in $\mathbb{D}$. Kayumov and Ponnusamy \cite{kayumov-2017} have obtained the following interesting result.
\begin{thm} \cite{kayumov-2017} \label{th-1.4}
	Let $f(z)=\sum_{n=0}^{\infty} a_{n}z^{n}$ be analytic in $\mathbb{D}$ and $|f(z)|\leq 1$. Then
	\begin{equation} \label{him-e-1.6}
		|f(z)|+\sum_{n=N}^{\infty}|a_n||z|^n\leq1
	\end{equation}
	for $|z|=r \leq R_{N}$, where $R_{N}$ is the positive root of the equation $\psi _{N}(r)=0$, $\psi _{N}(r)=2 (1+r)r^{N}-(1-r)^{2}$. The radius $R_{N}$ is the best possible. Moreover, 
	\begin{equation} \label{him-e-1.7}
		|f(z)|^{2}+\sum_{n=N}^{\infty}(|a_n|+|b_n|)|z|^n\leq1
	\end{equation}
	for $|z|=r\leq R'_{N}$, where $R'_{N}$ is the positive root of the equation $ (1+r)r^{N}-(1-r)^{2}$. The radius $R'_{N}$ is the best possible.
\end{thm}
In $ 2020 $, Ponnusamy {\it et al.} \cite{Ponnusamy-Vijaya-ResultsMath-2020} established the following refined Bohr inequality by applying a refined version of the coefficient inequalities.
\begin{thm} \cite{Ponnusamy-Vijaya-ResultsMath-2020} \label{th-1.6}
	Let $f(z)=\sum_{n=0}^{\infty} a_{n}z^{n}$ be analytic in $\mathbb{D}$ and $|f(z)|\leq 1$. Then
	$$
	\sum_{n=0}^{\infty} |a_{n}|r^{n}+ \left(\frac{1}{1+|a_{0}|}+\frac{r}{1-r}\right) \sum_{n=1}^{\infty}|a_n|^2r^{2n}\leq 1
	$$
	for $|z|=r\leq 1/(2+|a_{0}|)$ and the numbers $1/(1+|a_{0}|)$ and $1/(2+|a_{0}|)$ cannot be improved. Moreover, 
	$$
	|a_{0}|^{2}+ \sum_{n=1}^{\infty} |a_{n}|r^{n}+ \left(\frac{1}{1+|a_{0}|}+\frac{r}{1-r}\right) \sum_{n=1}^{\infty}|a_n|^2r^{2n}\leq 1
	$$
	for $|z|=r\leq 1/2$ and the numbers $1/(1+|a_{0}|)$ and $1/2$ cannot be improved.
\end{thm}
 A harmonic mapping in the unit disk $ \mathbb{D} $ is a complex-valued function $ f=u+iv $ of $ z=x+iy $ in $ \mathbb{D} $, which satisfies the Laplace equation $ \Delta f=4f_{z\bar{z}}=0 $, where $ f_{z}=(f_{x}-if_{y})/2 $ and $ f_{\bar{z}}=(f_{x}+if_{y})/2 $ and $ u $ and $ v $ are real-valued harmonic functions in $ \mathbb{D} $. It follows that the function $ f $ admits the canonical representation $ f=h+\bar{g} $, where $ h $ and $ g $ are analytic in $ \mathbb{D} $.  Let $ \mathcal{H} $ be the class of complex-valued harmonic mappings $ f=h+\bar{g}$ defined in the unit disk $ \mathbb{D} $, where $ h $ and $ g $ are analytic functions in $ \mathbb{D} $ with the normalization $ h(0)=0=h^{\prime}(0)-1 $ and $ g(0)=0 $. Let $ \mathcal{H}_{0}=\{f=h+\bar{g}\in\mathcal{H} : g^{\prime}(0)=0\}. $\vspace{1.5mm}

A domain $ \Omega\subseteq \mathbb{C} $ is called starlike with respect to a point $ z_0\in\Omega $ if the line segment joining $ z_0 $ to any other point in $ \Omega $ must lies in $ \Omega $. In particular, if $ z_0=0 $, then $ \Omega $ is called starlike. A complex-valued harmonic mapping $ f $ is said to be starlike if $ f(\mathbb{D}) $ is starlike. It is well-known that starlikeness is a heriditary property.  That is, if an anlytic function maps the unit disk $ \mathbb{D} $ univalently into a starlike domain, then it also maps the concentric circle onto a starlike domain. In general, starlike harmonic functions do not have this property. The failure of this hereditary property of starlike harmonic mappings led to the introduction of fully starlike harmonic mappings.\vspace{1.5mm}

It is well-known \cite{Nagpal-Ravi-APM-2013} that a sense-preserving harmonic function $ f=h+\overline{g} $ is fully starlike in $ \mathbb{D} $ if the analytic function $ h+\epsilon g $ is starlike in $ \mathbb{D} $ for each $ \epsilon $ with $ |\epsilon|=1 $.

\vspace{1.5mm}

Let 
$$\mathcal{P}^{0}_{\mathcal{H}}(M)=\{f=h+\overline{g} \in \mathcal{H}_{0}: \real (zh^{\prime\prime}(z))> -M+|zg^{\prime\prime}(z)|, \; z \in \mathbb{D}\; \mbox{and }\; M>0\},
$$ where \begin{equation}\label{e-1.2}
	f(z)=h(z)+\overline{g(z)}=z+\sum_{n=2}^{\infty}a_nz^n+\overline{\sum_{n=2}^{\infty}b_nz^n}.
\end{equation}

It is well-known that  \cite{Ghosh-Vasudevarao-BAMS-2020} functions in $ \mathcal{P}^{0}_{\mathcal{H}}(M) $  are starlike for $ 0<M<1/\log 4 $. In particular, functions in $ \mathcal{P}^{0}_{\mathcal{H}}(M) $ are fully starlike for $ 0<M<1/\log 4 $ (see \cite{Ghosh-Vasudevarao-BAMS-2020}). The growth theorem and sharp coefficient bounds for the functions in $ \mathcal{P}^{0}_{\mathcal{H}}(M) $ have been studied in \cite{Ghosh-Vasudevarao-BAMS-2020}. The following results on the coefficient bounds and growth estimation play vital role in proving our main results.
\begin{lem} \label{lem-1.8} \cite{Ghosh-Vasudevarao-BAMS-2020}
Let $f=h+\overline{g}\in \mathcal{P}^{0}_{\mathcal{H}}(M)$ for some $M>0$ and be of the form \eqref{e-1.2}. Then for $n\geq 2,$ 
\begin{enumerate}
\item[(i)] $\displaystyle |a_n| + |b_n|\leq \frac {2M}{n(n-1)}; $\\[2mm]

\item[(ii)] $\displaystyle ||a_n| - |b_n||\leq \frac {2M}{n(n-1)};$\\[2mm]

		\item[(iii)] $\displaystyle |a_n|\leq \frac {2M}{n(n-1)}.$
	\end{enumerate}
	The inequalities  are sharp with extremal function  $f_{M}$ given by $f_{M}(z)=z+ 2M \sum\limits_{n=2}^{\infty} \dfrac{z^n}{n(n-1)}.
	$
\end{lem}
\begin{lem}\cite{Ghosh-Vasudevarao-BAMS-2020}\label{lem-1.9}
	Let $f \in \mathcal{P}^{0}_{\mathcal{H}}(M)$. Then 
	\begin{equation} \label{e-1.4}
		|z| +2M \sum\limits_{n=2}^{\infty} \dfrac{(-1)^{n-1}|z|^{n}}{n(n-1)} \leq |f(z)| \leq |z| + 2M \sum\limits_{n=2}^{\infty} \dfrac{|z|^{n}}{n(n-1)}.
	\end{equation}
Both  inequalities are sharp for the function $f_{M}$ given by $f_{M}(z)=z+ 2M \sum\limits_{n=2}^{\infty} \dfrac{z^n}{n(n-1)}.
$
\end{lem}
 In $ 2017 $, Kayumov and Ponnusamy \cite{kayumov-2017} proved the following improved version of Bohr's inequality.
\begin{thm}\cite{kayumov-2017}\label{th-1.11}
	Let $f(z)=\sum_{n=0}^{\infty} a_{n}z^{n}$ be analytic in $\mathbb{D}$, $|f(z)|\leq 1$ and $S_{r}$ denotes the area of the image of the subdisk $|z|<r$ under mapping $f$. Then 
	\begin{equation}
		B_{1}(r):=\sum_{n=0}^{\infty}|a_n|r^n+ \frac{16}{9} \left(\frac{S_{r}}{\pi}\right) \leq 1 \quad \mbox{for} \quad r \leq \frac{1}{3}
	\end{equation}
	and the numbers $1/3$ and $16/9$ cannot be improved. Moreover, 
	\begin{equation}
		B_{2}(r):=|a_{0}|^{2}+\sum_{n=1}^{\infty}|a_n|r^n+ \frac{9}{8} \left(\frac{S_{r}}{\pi}\right) \leq 1 \quad \mbox{for} \quad r \leq \frac{1}{2}
	\end{equation}
	and the numbers $1/2$ and $9/8$ cannot be improved.
\end{thm}

The main aim of this paper is to obtain the sharp Bohr-Rogosinski inequality, improved Bohr inequality and Refined Bohr inequality and Bohr radius for the class $\mathcal{P}^{0}_{\mathcal{H}}(M)$. The organization of this paper as follows: In Section 2, we prove the sharp Bohr-Rogosinski, improved Bohr radius and Bohr radius for the class $ \mathcal{P}^{0}_{\mathcal{H}}(M) $. In Section 3, we give the proof of all the main results of this paper.
 \section{Main results}
 
  \par Before stating the main results of this paper, we recall the definiton of dilogarithm. The dilogarithm function $ {\rm Li}_2(z) $ is defined by the power series
  \begin{equation*}
  	{\rm Li}_2(z)=\sum_{n=1}^{\infty}\frac{z^n}{n^2}\;\; \mbox{for}\;\; |z|<1.
  \end{equation*}
  The definition and the name, of course, come from the analogy with the Taylor series of the ordinary logarithm around $ 1 $, 
  \begin{equation*}
  	-\log(1-z)=\sum_{n=1}^{\infty}\frac{z^n}{n}\;\;\mbox{for}\;\; |z|<1,
  \end{equation*}
  which leads similarly to the defintion of the polylogarithm
  \begin{equation*}
  	{\rm Li}_m(z)=\sum_{n=1}^{\infty}\frac{z^n}{n^m}\;\; \mbox{for}\;\; |z|<1,\;\; m=1, 2, 3, \ldots.
  \end{equation*}
  In addition, the relation 
  \begin{equation*}
  	\frac{d}{dz}\left({\rm Li}_m(z)\right)=\frac{1}{z}{\rm Li}_{m-1}(z)\;\; \mbox{for}\;\; m\geq 2
  \end{equation*}
  is obvious and leads by induction to the extension of the domain of definition of $ {\rm Li}_m $ to the cut plane $ \mathbb{C}\setminus [1, \infty) $. In particular, the analytic continuation of the dilogarithm is given by
  \begin{equation*}
  	{\rm Li}_2(z)=-\int_{0}^{z}\log (1-u)\frac{du}{u}\;\; \mbox{for}\;\; z\in\mathbb{C}\setminus [1,\infty).
  \end{equation*}
Using Lemma \ref{lem-1.8} and Lemma \ref{lem-1.9}, for the functions in the class $ \mathcal{P}^{0}_{\mathcal{H}}(M) $, we prove the following sharp Bohr inequality which is an harmonic analogue of Theorem \ref{th-1.4}.
 \begin{thm}\label{th-2.1}
 	Let $ f\in \mathcal{P}^{0}_{\mathcal{H}}(M) $ be given \eqref{e-1.2} with $ 0<M<1/(2(\ln 4-1)) $. Then 
 	\begin{equation*}
 		|z|+|f(z)|+\sum_{n=2}^{\infty}(|a_n|+|b_n|)|z|^n\leq d\left(f(0),\partial f(\mathbb{D})\right)
 	\end{equation*} 
 	for $ |z|=r\leq r_{_M} $, where $ r_{_M} $ is the unique root of the equation
 	\begin{equation}\label{e-2.3}
 	2r+4M\left(r+(1-r)\log(1-r)\right)-1-2M(1-2\log 2)=0
 	\end{equation}
 in $ (0,1) $. The radius $ r_{_M} $ is the best possible.
 \end{thm}
\begin{rem}
	For particular values of $ M $ in Theorem \ref{th-2.1}, we have $$ 0< M<{1}/{(2(2\ln 2-1))}\approx 1.294349.$$ A simple computation gives  the Bohr radii $ r_{0.1}\approx 0.438485 $, $ r_{0.2}\approx 0.387786 $, $ r_{0.3}\approx 0.343722 $, $ r_{0.4}\approx 0.304054 $, $ r_{0.5}\approx 0.267404 $, $ r_{0.6}\approx 0.23283 $, $ r_{0.7}\approx 0.19963 $, $ r_{0.8}\approx 0.167229 $, $ r_{0.9}\approx 0.135111 $, $ r_{1.0}\approx 0.102764 $, $ r_{1.25}\approx 0.0167782 $.
\end{rem}
 \begin{figure}[!htb]
 	\begin{center}
 		\includegraphics[width=0.65\linewidth]{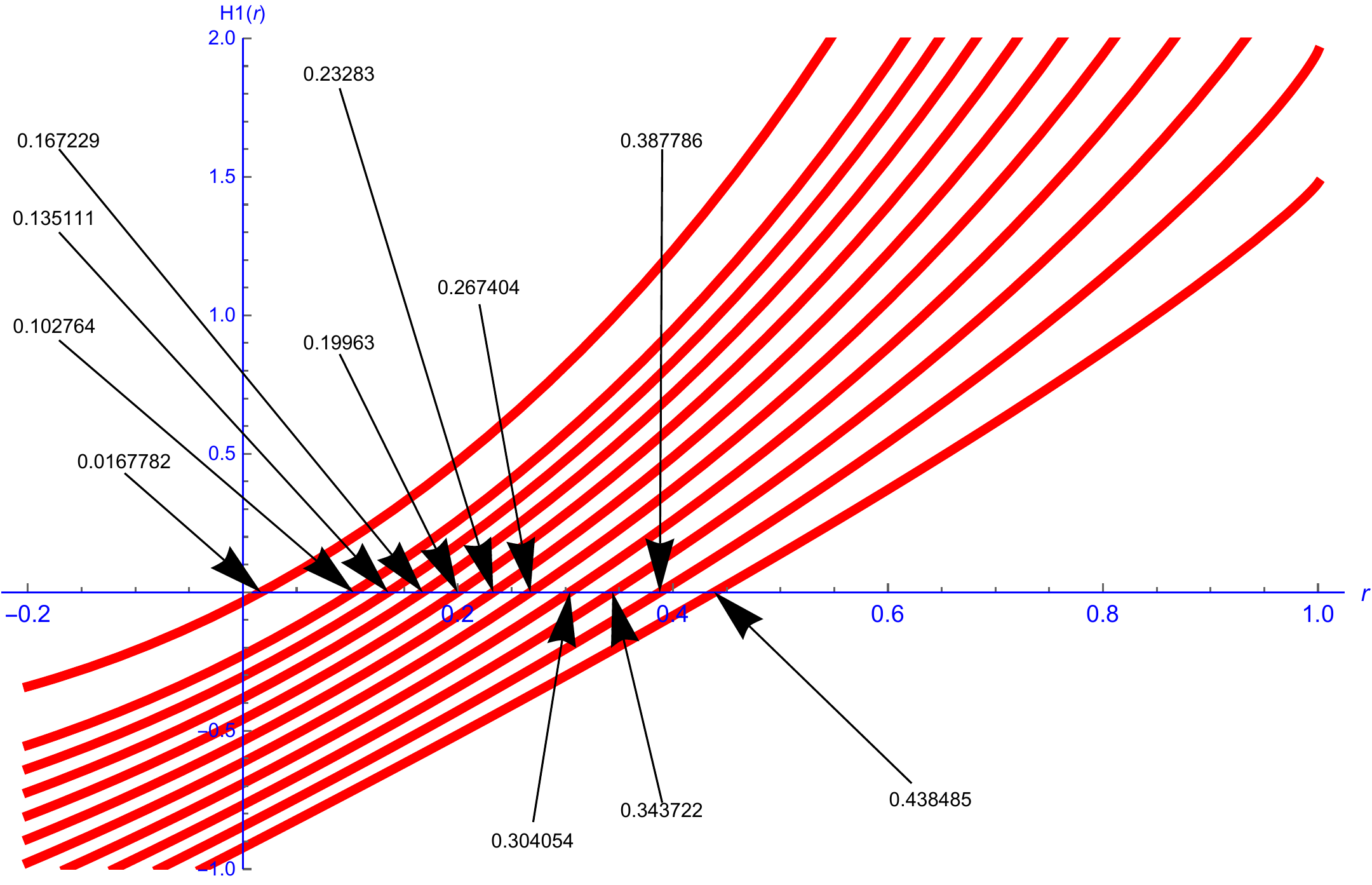}
 	\end{center}
 	\caption{The roots of \eqref{e-2.3} for different values of $ M $ when $0< M<{1}/{(2(2\ln 2-1))} $.}
 \end{figure}
Figure 1 illustrates the roots of \eqref{e-2.3} for different values of $ M $ when $0< M<{1}/{(2(2\ln 2-1))} $. Let $ S_r $ denote the area of the image of of the subdisk $ \mathbb{D}_r:=\{z\in\mathbb{C} : |z|<r\} $ under the mappings in the class $ \mathcal{P}^{0}_{\mathcal{H}}(M) $. We now prove the sharp Bohr inequality which is a harmonic analogue of Theorem \ref{th-1.11}.
 \begin{thm}\label{th-2.3}
Let $ f\in \mathcal{P}^{0}_{\mathcal{H}}(M) $ be given \eqref{e-1.2} with $ 0<M<1/(2(\ln 4-1)) $. Then 
\begin{equation*}
	|f(z)|+\sum_{n=2}^{\infty}(|a_n|+|b_n|)|z|^n+\frac{S_r}{\pi}\leq d\left(f(0),\partial f(\mathbb{D})\right)
\end{equation*} 
for $ |z|=r\leq r^*_{_M} $, where $ r^*_{_M} $ is the unique root of the equation
\begin{align}
	\label{e-2.4}
	&r+4M\left(r+(1-r)\log(1-r)\right)+4M^2\left(r^2{\rm Li}_2(r^2)-r^2-\left(1-r^2\right)\log\left(1-r^2\right)\right)\\&\nonumber\quad\quad-1-2M(1-2\log 2)=0
\end{align}
in $ (0,1) $. The radius $ r^*_{_M} $ is the best possible.
\end{thm}
\begin{rem}
	For particular values of $ M $ in Theorem \ref{th-2.3}, a simple computation yields the Bohr radii $ r^*_{0.1}\approx 0.546723 $, $ r^*_{0.2}\approx 0.487374 $, $ r^*_{0.3}\approx 0.435926 $, $ r^*_{0.4}\approx 0.389886 $, $ r^*_{0.5}\approx 0.0.347564 $, $ r^*_{0.6}\approx 0.307711 $, $ r^*_{0.7}\approx 0.269313 $, $ r^*_{0.8}\approx 0.231445 $, $ r^*_{0.9}\approx 0.193148, $ $ r^*_{1.0}\approx 0.153247 $, $ r^*_{1.25}\approx 0.0371406 $.
\end{rem}
\begin{figure}[!htb]
	\begin{center}
		\includegraphics[width=0.65\linewidth]{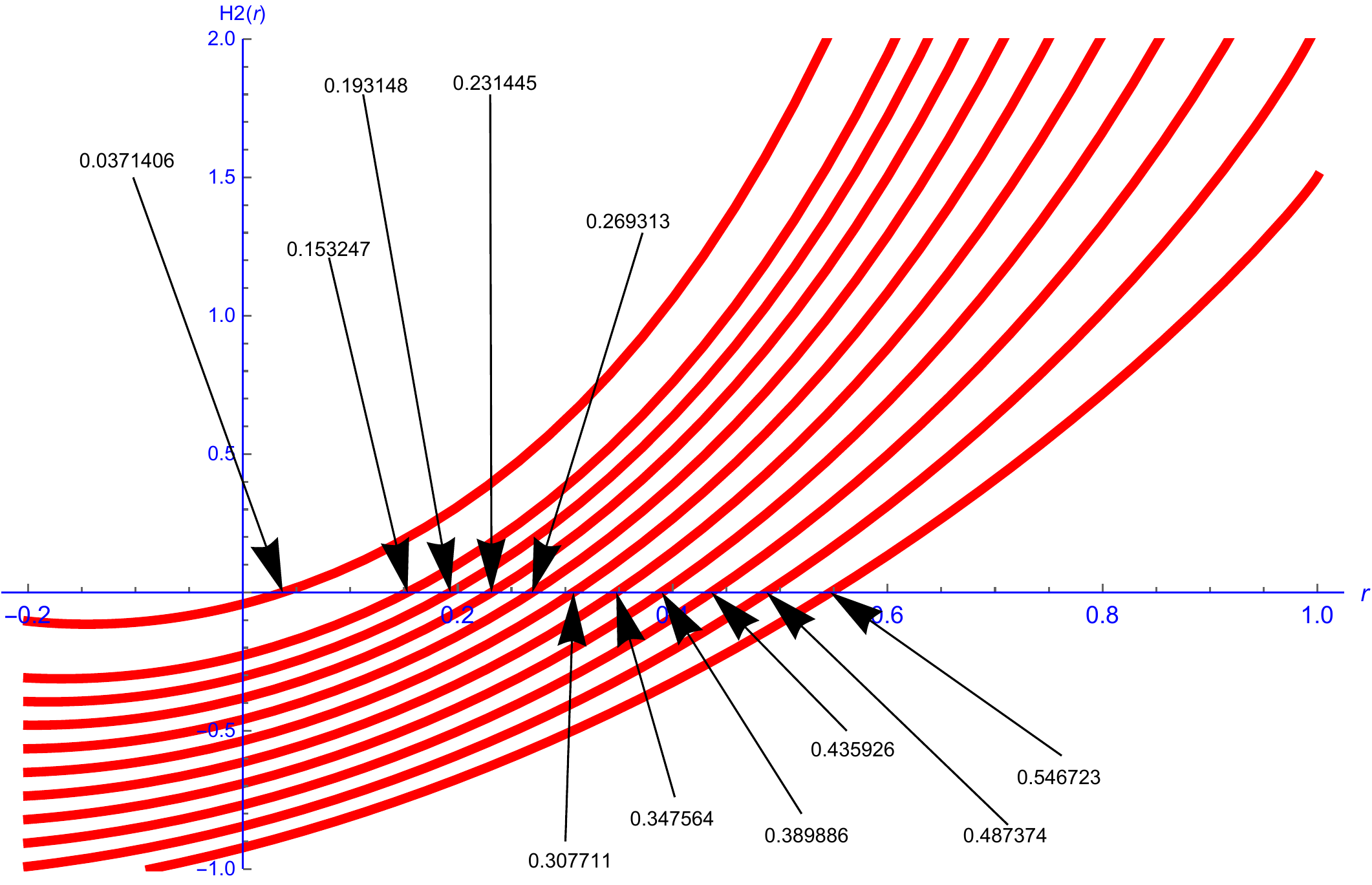}
	\end{center}
	\caption{The roots of \eqref{e-2.4} for different values of $ M $ when $0< M<{1}/{(2(2\ln 2-1))} $.}
\end{figure} 
Figure 2 illustrates the roots of \eqref{e-2.4} for different values of $ M $ when $0< M<{1}/{(2(2\ln 2-1))} $. The Jacobian of a complex-valued functions of the form $ f=h+\bar{g} $ is defined by $ J_f(z)=|h^{\prime}(z)|^2-|g^{\prime}(z)|^2 $. It is known that $ f $ is sense-preserving if $ J_f(z)>0 $ and $ f $ is sense-reversing if $ J_f(z)<0 $ in $ \mathbb{D} $. Using Lemma \ref{lem-1.8} and Lemma \ref{lem-1.9}, we prove the following result for the functions in the class $ \mathcal{P}^{0}_{\mathcal{H}}(M) $.
\begin{thm}\label{th-2.5}
	Let $ f\in \mathcal{P}^{0}_{\mathcal{H}}(M) $ be given \eqref{e-1.2} with $ 0<M<1/(2(\ln 4-1)) $. Then 
	\begin{equation*}
		\sum_{n=2}^{\infty}(|a_n|+|b_n|)|z|^n+\sqrt{|J_f(z)|}|z|\leq d\left(f(0),\partial f(\mathbb{D})\right)
	\end{equation*} 
	for $ |z|=r\leq r^{**}_{_M} $, where $ r^{**}_{_M} $ is the unique root of the equation
	\begin{align}
		\label{e-2.5}
		&2r+2M\left(r+(1-r)\log(1-r)\right)-2Mr\log(1-r)-1-2M(1-2\log 2)=0
	\end{align}
	in $ (0,1) $. The radius $ r^{**}_{_M} $ is the best possible.
\end{thm}
\begin{figure}[!htb]
	\begin{center}
		\includegraphics[width=0.65\linewidth]{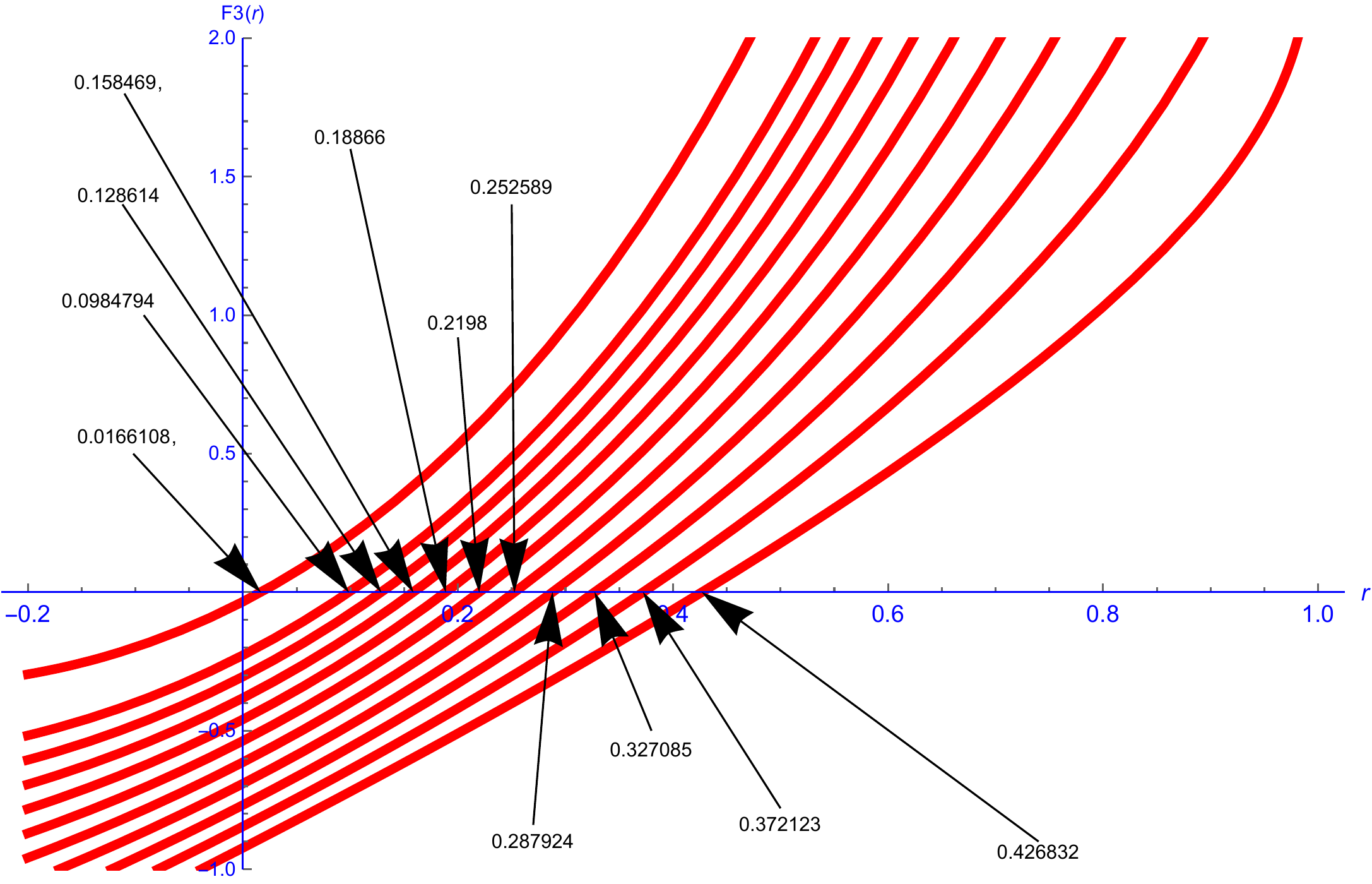}
	\end{center}
	\caption{The roots of \eqref{e-2.5} for different values of $ M $ when $0< M<{1}/{(2(2\ln 2-1))} $.}
\end{figure} 
 \begin{rem}
 	For particular values of $ M $ in Theorem \ref{th-2.5}, a simple computation yields the Bohr radii $ r^{**}_{0.1}\approx 0.426832 $, $ r^{**}_{0.2}\approx 0.372123 $, $ r^{**}_{0.3}\approx 0.327085 $, $ r^{**}_{0.4}\approx 0.287924 $, $ r^{**}_{0.5}\approx 0.252589 $, $ r^{**}_{0.6}\approx 0.2198 $, $ r^{**}_{0.7}\approx 0.18866 $, $ r^{**}_{0.8}\approx 0.158469 $, $ r^{**}_{0.9}\approx 0.128614 $, $ r^{**}_{1.0}\approx 0.0984794 $, $ r^{**}_{1.25}\approx 0.0166108 $.
 \end{rem}
Figure 3 illustrates the roots of \eqref{e-2.5} for different values of $ M $ when $0< M<{1}/{(2(2\ln 2-1))} $. We prove the following refine sharp bohr inequality for functions in the class $ \mathcal{P}^{0}_{\mathcal{H}}(M) $ which is a harmonic analogue version of Theorem \ref{th-1.6} .
\begin{thm}\label{th-2.7}
	Let $ f\in \mathcal{P}^{0}_{\mathcal{H}}(M) $ be given \eqref{e-1.2} with $ 0<M<1/(2(\ln 4-1)) $. Then 
	\begin{equation*}
		|f(z)|^2+\sum_{n=2}^{\infty}(|a_n|+|b_n|)|z|^n+\frac{r}{1-r}\sum_{n=2}^{\infty}(|a_n|+|b_n|)^2|z|^{2n}\leq d\left(f(0),\partial f(\mathbb{D})\right)
	\end{equation*} 
	for $ |z|=r\leq r^{\prime}_{_M} $, where $ r^{\prime}_{_M} $ is the unique root of the equation
	\begin{align}
		\label{e-2.6}
		&\left(r+2M\left(r+(1-r)\log(1-r)\right)\right)^2+2M\left(r+(1-r)\log(1-r)\right)\\&\nonumber\quad\quad-2Mr\log(1-r)-1-2M(1-2\log 2)=0
	\end{align}
	in $ (0,1) $. The radius $ r^{\prime}_{_M} $ is the best possible.
\end{thm}
\begin{figure}[!htb]
	\begin{center}
		\includegraphics[width=0.65\linewidth]{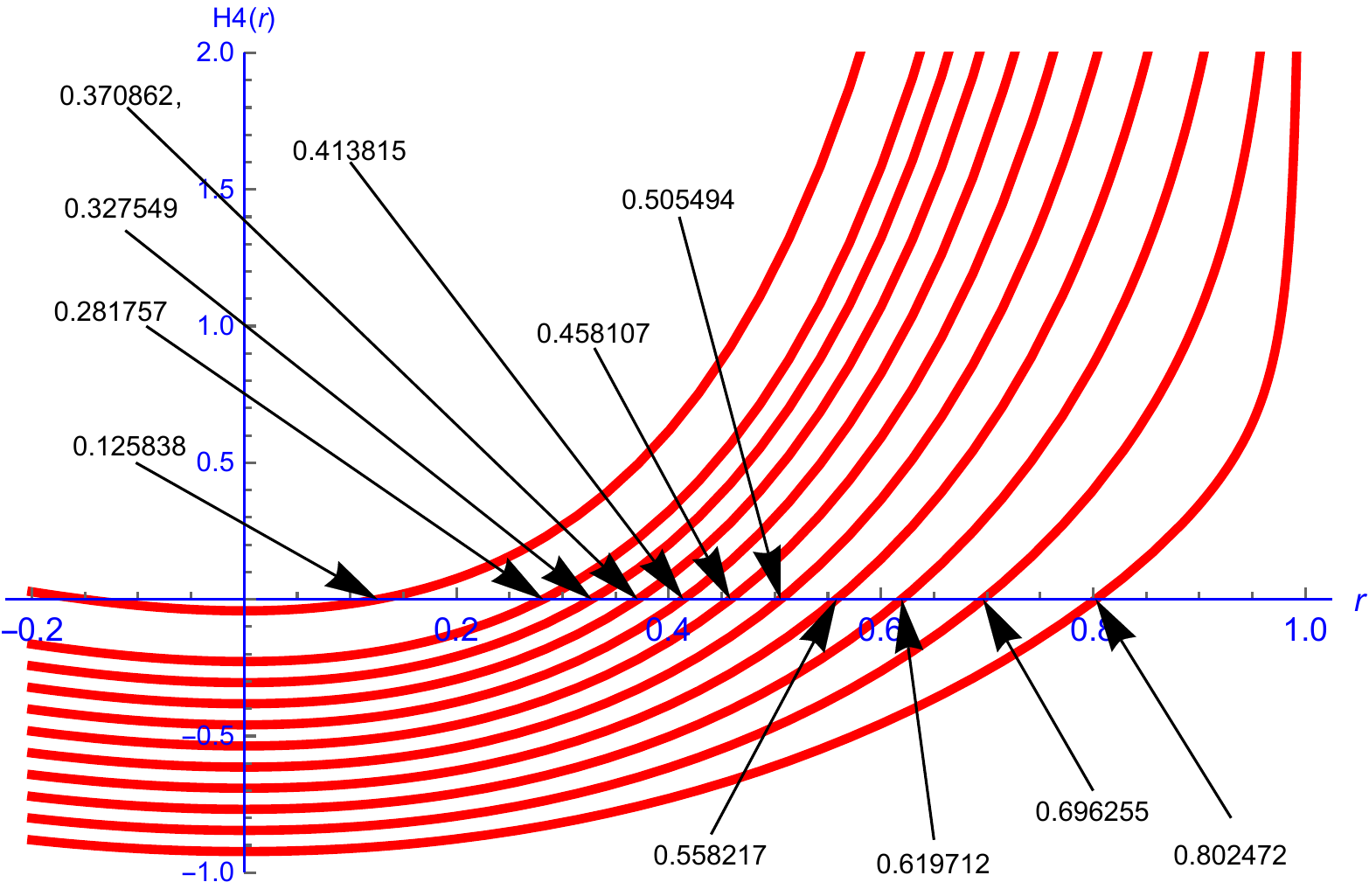}
	\end{center}
	\caption{The roots of \eqref{e-2.6} for different values of $ M $ when $0< M<{1}/{(2(2\ln 2-1))} $.}
\end{figure} 
\begin{rem}
	For particular values of $ M $ in Theorem \ref{th-2.7}, a simple computation gives the Bohr radii as follow: $ r^{\prime}_{0.1}\approx 0.802472 $, $ r^{\prime}_{0.2}\approx 0.696255 $, $ r^{\prime}_{0.3}\approx 0.619712 $, $ r^{\prime}_{0.4}\approx 0.558217 $, $ r^{\prime}_{0.5}\approx 0.505494 $, $ r^{\prime}_{0.6}\approx 0.458107 $, $ r^{\prime}_{0.7}\approx 0.413815 $, $ r^{\prime}_{0.8}\approx 0.370862 $, $ r^{\prime}_{0.9}\approx 0.327549 $, $ r^{\prime}_{1.0}\approx 0.281757 $, $ r^{\prime}_{1.25}\approx 0.125838 $.
\end{rem}
Figure 4 illustrates the roots of \eqref{e-2.6} for different values of $ M $ when $0< M<{1}/{(2(2\ln 2-1))} $
\section{Proof of the main results}	
 For $ f \in \mathcal{P}^{0}_{\mathcal{H}}(M) $, we have
\begin{equation}\label{e-3.1}
	|f(z)|\geq |z|+2M\sum\limits_{n=2}^{\infty}  \dfrac{(-1)^{n-1}}{n(n-1)}|z|^{n} \quad \mbox{for } \quad |z|<1.
\end{equation}
Then the Euclidean distance between $f(0)$ and the boundary of $f(\mathbb{D})$ is
\begin{equation}\label{e-3.2}
	d(f(0), \partial f(\mathbb{D}))= \liminf\limits_{|z|\rightarrow 1} |f(z)-f(0)|.
\end{equation}
Therefore, a simple computation using \eqref{e-1.4} and \eqref{e-3.2} shows that 
\begin{equation}\label{e-3.3}
	d(f(0), \partial f(\mathbb{D})) \geq 1+2M\sum\limits_{n=2}^{\infty}  \dfrac{(-1)^{n-1}}{n(n-1)}=1+2M(1-\ln 4).
\end{equation}
\begin{pf}[\bf Proof of Theorem \ref{th-2.1}]
	Let $f \in \mathcal{P}^{0}_{\mathcal{H}}(M)$ be given by \eqref{e-1.2}. Using Lemma \ref{lem-1.8} and Lemma \ref{lem-1.9}, for $ |z|=r $, we obtain 
	\begin{align}\label{e-3.4}
		|z|+|f(z)|+\sum_{n=2}^{\infty}(|a_n|+|b_n|)|z|^n &\leq 2r+\sum_{n=2}^{\infty}\frac{2Mr^n}{n(n-1)}+\sum_{n=2}^{\infty}\frac{2Mr^n}{n(n-1)}\\&=2r+4M\sum_{n=2}^{\infty}\frac{r^n}{n(n-1)}\nonumber.
	\end{align}
A simple computation shows that
\begin{align}\label{e-3.5}
	\sum_{n=2}^{\infty}\frac{r^n}{n(n-1)}=r\sum_{n=1}^{\infty}\frac{r^n}{n}+\sum_{n=2}^{\infty}-\frac{r^n}{n}= r+(1-r)\log(1-r)
\end{align}
and 
\begin{equation*}
	\sum_{n=2}^{\infty}\frac{(-1)^{n-1}}{n(n-1)}=1-2\log 2\approx -0.386294.
\end{equation*}
Therefore, using \eqref{e-3.5} in \eqref{e-3.4}, we obtain
\begin{align*}
	|z|+|f(z)|+\sum_{n=2}^{\infty}(|a_n|+|b_n|)r^n\leq 2r+4M\left(r+(1-r)\log(1-r)\right). 
\end{align*}
It is easy to see that
	\begin{equation*}
	2r+4M\left(r+(1-r)\log(1-r)\right)\leq 1+2M(1-2\ln 2)
	\end{equation*}
	for $ r\leq r_{_M} $, where $ r_{_M} $ is a root of $ H_1(r)=0 $ in $ (0,1) $, where
$ H_1 : [0,1)\rightarrow \mathbb{R} $ is defined by 
\begin{align*}
		H_1(r):&=2r+4M\left(r+(1-r)\log(1-r)\right)-1-2M(1-2\ln 2)\\&=2r+4M\sum_{n=2}^{\infty}\frac{r^n}{n(n-1)}-1-2M	\sum_{n=2}^{\infty}\frac{(-1)^{n-1}}{n(n-1)}.
\end{align*}
The existence of the root $ r_{_M} $ is guaranteed by the Intermediate value theorem since $ H_1 $ is a continuous function on $ [0,1] $, differentiable on $ (0,1) $ and $ H_1(0)=-1-2M(1-2\ln 2)<0 $ for $0< M<1/(2(\ln 4-1)) $ and $ \displaystyle\lim_{r\rightarrow 1}H_1(r)=+\infty. $ To show that $ r_{_M} $ is the unique root of $ H_1(r)=0 $, it is enough to show that $ H_1 $ is strictly monotone function on $ (0,1) $. A simple computation shows that
\begin{equation*}
	\frac{d}{dr}\left(H_1(r)\right)=2+4M\sum_{n=2}^{\infty}\frac{r^{n-1}}{n-1}>0
\end{equation*}
for $ r\in (0,1) $, which shows that $ H_1 $ is strictly increasing in $ (0,1) $. Hence, the root $ r_{_M} $ is the unique root of $ H_1 $ in $ (0,1) $ and hence, we have
\begin{equation}\label{e-3.6}
2r_{_M}+4M\left(r_{_M}+(1-r_{_M})\log(1-r_{_M})\right)=1+2M(1-2\ln 2).
\end{equation}
To show that $ r_{_M} $ is the best possible, we consider the function $ f=f_M $ is given by 
	\begin{equation}\label{e-3.7}
		f_{M}(z)=z+\sum_{n=2}^{\infty}\frac{2Mz^n}{n(n-1)}.
	\end{equation}
	It is easy to see that $ f_{M}\in\mathcal{P}^{0}_{\mathcal{H}}(M) $ and for $ f=f_{M} $, we have 
	\begin{equation}\label{e-3.8}
		d(f_M(0),\partial f_M(\mathbb{D}))=1+2M(1-2\ln 2).
	\end{equation}
	For the function $ f=f_{M} $ and $ z=r $, a simple computation shows that 
	\begin{align} \label{e-3.9}
		r+|f_M(r)|+\sum_{n=2}^{\infty}(|a_n|+|b_n|)r^n=2r+\sum_{n=2}^{\infty}\frac{4Mr^n}{n(n-1)}.
	\end{align}
Therefore, a simple computation using \eqref{e-3.7}, \eqref{e-3.8} and \eqref{e-3.9} for $ f=f_M $ and $ r>r_{_M} $ shows that
\begin{align*}
	&|z|+|f_M(z)|+\sum_{n=2}^{\infty}(|a_n|+|b_n|)|z|^n\\ &= 2r+\sum_{n=2}^{\infty}\frac{2Mr^n}{n(n-1)}+\sum_{n=2}^{\infty}\frac{2Mr^n}{n(n-1)}\\&>2r_{_M}+\sum_{n=2}^{\infty}\frac{2Mr^n_{_M}}{n(n-1)}+\sum_{n=2}^{\infty}\frac{2Mr^n_{_M}}{n(n-1)}\\&=2r+4M\sum_{n=2}^{\infty}\frac{r^n_{_M}}{n(n-1)}\\&=2r_{_M}+4M\left(r_{_M}+(1-r_{_M})\log(1-r_{_M})\right)\\&=1+2M(1-2\ln 2)\\&=d(f_M(0),\partial f_M(\mathbb{D}))
\end{align*}
and hence $ r_{_M} $ is the best possible. This completes the proof.
\end{pf}
\begin{pf}[\bf Proof of Theorem \ref{th-2.3}]
	Let $f \in \mathcal{P}^{0}_{\mathcal{H}}(M)$ be given by \eqref{e-1.2}. It is well-known that (see \cite[p.113]{Duren-2004}) the area of the disk $ \mathbb{D}_r:=\{z\in\mathbb{C} : |z|<r\} $ under the harmonic map $ f=h+\bar{g} $ is 
	\begin{align}\label{e-3.18}
		S_r=\iint\limits_{\mathbb{D}_r}J_f(z)dxdy=\iint\limits_{\mathbb{D}_r}\left(|{h^{\prime}(z)}|^2-|{g^{\prime}(z)}|^2\right)dxdy.
	\end{align}
It is easy to see that
	\begin{align}\label{e-3.19}
		\iint\limits_{\mathbb{D}_r}|h^{\prime}(z)|^2dxdy&=\int_{0}^{r}\int_{0}^{2\pi}|h^{\prime}(\rho e^{i\theta})|^2\rho d\theta d\rho\\&=\nonumber\int_{0}^{r}\int_{0}^{2\pi}\rho\left(\sum_{n=1}^{\infty}n\;a_n\rho^{n-1}e^{i(n-1)\theta}\right)\left(\sum_{n=1}^{\infty}n\;\bar{a}_n\rho^{n-1}e^{-i(n-1)\theta}\right)d\theta d\rho\\&=\nonumber \int_{0}^{r}\left(\sum_{n=1}^{\infty}2\pi n^2 |a_n|^2\rho^{2n-1}\right)d\rho\\&\nonumber=\sum_{n=1}^{\infty}2\pi n^2|a_n|^2\frac{r^{2n}}{2n}\\&\nonumber=\pi \sum_{n=1}^{\infty}n|a_n|^2r^{2n}.
	\end{align}
	Similarly, for $ g(z)=\sum_{n=2}^{\infty}b_nz^n $, we obtain
	\begin{align}\label{e-3.20}
		\iint\limits_{\mathbb{D}_r}&|g^{\prime}(z)|^2dxdy=\pi\sum_{n=2}^{\infty}n|b_n|^2r^{2n}.
	\end{align}
	In view of Lemma \ref{lem-1.8}, \eqref{e-3.18}, \eqref{e-3.19} and  \eqref{e-3.20}, we obtain
	\begin{align}\label{e-3.21}
		\frac{S_r}{\pi}&=\frac{1}{\pi}\iint\limits_{\mathbb{D}_r}\left(|h^{\prime}(z)|^2-|g^{\prime}(z)|^2\right)dxdy\\&\nonumber=r^2+\sum_{n=2}^{\infty}n\left(|a_n|^2-|b_n|^2\right)r^{2n}\\&\nonumber=r^2+\sum_{n=2}^{\infty}n(|a_n|+|b_n|)(|a_n|-|b_n|)r^{2n}\\&\nonumber=r^2+\sum_{n=2}^{\infty}\frac{4M^2r^{2n}}{n\left(n-1\right)^2}.
	\end{align}
Therefore, in view of Lemma \ref{lem-1.8}, Lemma \ref{lem-1.9} and \eqref{e-3.20}, for $ |z|=r $ yields
\begin{align}\label{e-3.10}
	&|f(z)|+\sum_{n=2}^{\infty}\left(|a_n|+|b_n|\right)|z|^n+\frac{S_r}{\pi}\\&\nonumber\leq r+4M\sum_{n=2}^{\infty}\frac{r^n}{n(n-1)}+\left(r^2+\sum_{n=2}^{\infty}\frac{4M^2r^{2n}}{n(n-1)^2}\right).
\end{align}
 A simple computation shows that
 \begin{align}\label{e-3.11}
 	\sum_{n=2}^{\infty}\frac{r^{2n}}{n(n-1)^2}&=\sum_{n=2}^{\infty}\frac{r^{2n}}{n}+\sum_{n=2}^{\infty}\frac{r^{2n}}{(n-1)^2}+\sum_{n=2}^{\infty}-\frac{r^{2n}}{n-1}\\&\nonumber=-r^2+\sum_{n=1}^{\infty}\frac{r^{2n}}{n}+r^2\sum_{n=2}^{\infty}\frac{r^{2n-2}}{(n-1)^2}+r^2\sum_{n=2}^{\infty}-\frac{r^{2n-2}}{n-1}\\&\nonumber=-r^2-\log(1-r^2)+r^2\sum_{n=1}^{\infty}\frac{r^{2n}}{n^2}-r^2\sum_{n=1}^{\infty}\frac{r^{2n}}{n}\\&=r^2{\rm Li}_2\left(r^2\right)-r^2-\left(1-r^2\right)\log\left(1-r^2\right).\nonumber
 \end{align}
Therefore, using \eqref{e-3.5} and \eqref{e-3.11}, from \eqref{e-3.10}, we obtain
\begin{align*}
	&|f(z)|+\sum_{n=2}^{\infty}\left(|a_n|+|b_n|\right)|z|^n+\frac{S_r}{\pi}\\&\leq r^2+r+4M\left(r+(1-r)\log(1-r)\right)+4M^2\left(r^2{\rm Li}_2\left(r^2\right)-r^2-\left(1-r^2\right)\log\left(1-r^2\right)\right).
\end{align*}
It is easy to see that
\begin{align*}
	&r^2+r+4M\left(r+(1-r)\log(1-r)\right)+4M^2\left(r^2{\rm Li}_2\left(r^2\right)-r^2-\left(1-r^2\right)\log\left(1-r^2\right)\right)\\&\leq 1+2M\left(1-2\log 2\right)
\end{align*}
for $ r\leq r^*_{_M} $, where $ r^*_{_M} $ is a root of $ H_2(r)=0 $, where $ H_2 : [0,1]\rightarrow\mathbb{R} $ is defined by
\begin{align*}
	H_2(r) :&=r^2+r+4M\left(r+(1-r)\log(1-r)\right)-1-2M\left(1-2\log 2\right)\\&\quad\quad+4M^2\left(r^2{\rm Li}_2\left(r^2\right)-r^2-\left(1-r^2\right)\log\left(1-r^2\right)\right)\\&=r^2+r+4M\sum_{n=2}^{\infty}\frac{r^n}{n(n-1)}+\sum_{n=2}^{\infty}\frac{4M^2r^{2n}}{n(n-1)^2}-1-2M\sum_{n=2}^{\infty}\frac{(-1)^{n-1}}{n(n-1)}.
\end{align*}
A simple computation shows that $ H_2(0)=-1-2M\left(1-2\log 2\right)<0 $
for $ 0<M<1/(2(\ln 4-1)) $ and 
\begin{align*}
	H_2(1)&=2+4M\sum_{n=2}^{\infty}\frac{1}{n(n-1)}+4M^2\sum_{n=2}^{\infty}\frac{1}{n(n-1)^2}-2M(1-2\log 2)\\&=2+4M\sum_{n=2}^{\infty}\left(\frac{1}{n-1}-\frac{1}{n}\right)+4M^2\left(\sum_{n=2}^{\infty}\frac{1}{(n-1)^2}+\sum_{n=2}^{\infty}\left(\frac{1}{n}-\frac{1}{n-1}\right)\right)\\&\quad\quad-2M(1-2\log 2)\\&=2+4M\left(\left(1-\frac{1}{2}\right)+\left(\frac{1}{2}-\frac{1}{3}\right)+\left(\frac{1}{3}-\frac{1}{4}\right)+\cdots\right)-2M(1-2\log 2)\\&\quad\quad+4M^2\left(\left(\frac{1}{2}-1\right)+\left(\frac{1}{3}-\frac{1}{2}\right)+\left(\frac{1}{4}\right)-\frac{1}{3}+\cdots+\sum_{n=1}^{\infty}\frac{1}{n^2}\right)\\&= 2+2M+4M^2\left(\frac{\pi^2}{6}-1\right)+4M\log 2\\&>0.
\end{align*}

 It is easy to see that $ H_2 $ is strictly increasing in $ (0,1) $ because
\begin{align*}
	\frac{d}{dr}\left(H_2(r)\right)=2r+1+4M\sum_{n=2}^{\infty}\frac{r^{n-1}}{n-1}+8M^2\sum_{n=2}^{\infty}\frac{r^{2n}}{(n-1)^2}>0
\end{align*}
for $ r\in (0,1) $. By the similar argument that is being used in the proof of Theorem \ref{th-2.1}, we can show that $ r^*_{_M} $ is the unique root of  $ H_2 $ in $ (0,1) $ and hence, we have \begin{align}\label{e-3.12}
		&\left(r^*_{_M}\right)^2+r^*_{_M}+4M\left(r^*_{_M}+(1-r^*_{_M})\log(1-r)\right)\\\nonumber&\quad+4M^2\left(\left(r^*_{_M}\right)^2{\rm Li}_2\left(\left(r^*_{_M}\right)^2\right)-\left(r^*_{_M}\right)^2-\left(1-\left(r^*_{_M}\right)^2\right)\log\left(1-\left(r^*_{_M}\right)^2\right)\right)\\&\nonumber=1+2M\left(1-2\log 2\right).
\end{align}
To show that $  r^*_{_M} $ is the best possible radius, we consider $ f=f_M $ given by \eqref{e-3.7}. A simple computation, using \eqref{e-3.8}, \eqref{e-3.10} and \eqref{e-3.12} for $ f=f_M $ and $ r^*_{_M}<r $ shows that
\begin{align*}
	&|f(z)|+\sum_{n=2}^{\infty}\left(|a_n|+|b_n|\right)|z|^n+\frac{S_r}{\pi}\\&\nonumber= r^2+r+4M\sum_{n=2}^{\infty}\frac{r^n}{n(n-1)}+\left(r^2+\sum_{n=2}^{\infty}\frac{4M^2r^{2n}}{n(n-1)^2}\right)\\&>\left(r^*_{_M}\right)^2+r^*_{_M}+4M\sum_{n=2}^{\infty}\frac{\left(r^*_{_M}\right)^n}{n(n-1)}+\left(\left(r^*_{_M}\right)^2+\sum_{n=2}^{\infty}\frac{4M^2\left(r^*_{_M}\right)^{2n}}{n(n-1)^2}\right)\\&=\left(r^*_{_M}\right)^2+r^*_{_M}+4M\left(r^*_{_M}+(1-r^*_{_M})\log(1-r)\right)\\\nonumber&\quad+4M^2\left(\left(r^*_{_M}\right)^2{\rm Li}_2\left(\left(r^*_{_M}\right)^2\right)-\left(r^*_{_M}\right)^2-\left(1-\left(r^*_{_M}\right)^2\right)\log\left(1-\left(r^*_{_M}\right)^2\right)\right)\\&\nonumber=1+2M\left(1-2\log 2\right)\\&=d(f_M(0),\partial f_M(\mathbb{D}))
\end{align*}
and hence  $ r^*_{_M} $ is the best possible radius. This completes the proof.
\end{pf}
\begin{pf}[\bf Proof of Theorem \ref{th-2.5}]
	Let $f \in \mathcal{P}^{0}_{\mathcal{H}}(M)$ be given by \eqref{e-1.2}. The Jacobian of a complex-valued function $ f \in \mathcal{P}^{0}_{\mathcal{H}}(M)$ has the following property
	\begin{align}\label{e-3.13}
	|J_f(z)|=|h^{\prime}(z)|^2-|g^{\prime}(z)|^2\leq |h^{\prime}(z)|^2\leq \left(1+2M\sum_{n=1}^{\infty}\frac{r^n}{n}\right)^2.
	\end{align}
	In view of Lemma \ref{lem-1.8}, Lemma \ref{lem-1.9} and \eqref{e-3.13}, for $ |z|=r $, we obtain 
	\begin{align}\label{e-3.14}
	&\sum_{n=2}^{\infty}\left(|a_n|+|b_n|\right)|z|^n+\sqrt{|J_f(z)|}|z|\\&\nonumber\leq r+\sum_{n=2}^{\infty}\frac{2Mr^n}{n(n-1)}+\left(1+2M\sum_{n=1}^{\infty}\frac{r^n}{n}\right)r.
	\end{align}
A simple computation shows that
\begin{equation}\label{e-3.15}
\begin{cases}
	\displaystyle\sum_{n=2}^{\infty}\frac{r^n}{n(n-1)}=r+(1-r)\log(1-r),\vspace{1.5mm}\\ \displaystyle\sum_{n=1}^{\infty}\frac{r^n}{n}=-\log(1-r).
\end{cases}
\end{equation}
Therefore, using \eqref{e-3.15} in \eqref{e-3.14}, we obtain
\begin{align*}
	\sum_{n=2}^{\infty}\left(|a_n|+|b_n|\right)|z|^n+\sqrt{|J_f(z)}|z|\leq 2r+2M\left(r+(1-2r)\log(1-r)\right).
\end{align*}
It is easy to see that
\begin{align*}
	2r+2M\left(r+(1-2r)\log(1-r)\right)\leq 1+2M(1-2\log 2)
\end{align*}
for $ r\leq r^{**}_{_M} $, where $ r^{**}_{_M} $ is a root of $ H_3(r)=0 $, where $ H_3 : [0,1]\rightarrow\mathbb{R} $ is defined by 
\begin{align*}
	H_3(r):&=2r+2M\left(r+(1-2r)\log(1-r)\right)-1-2M(1-2\log 2)\\&=r+\sum_{n=2}^{\infty}\frac{2Mr^n}{n(n-1)}+\left(1+2M\sum_{n=1}^{\infty}\frac{r^n}{n}\right)r-1-2M\sum_{n=2}^{\infty}\frac{(-1)^{n-1}}{n(n-1)}.
\end{align*}
With the help of the similar argument that is being used in the proof of Theorem \ref{th-2.1} and Theorem \ref{th-2.3}, we can show that $ r^{**}_{_M} $ is the unique root of $ H_3(r) $ in $ (0,1) $. Thus, we have
\begin{equation}\label{e-3.16}
	2r^{**}_{_M}+2M\left(r^{**}_{_M}+(1-2r^{**}_{_M})\log(1-r)\right)=1+2M(1-2\log2).
\end{equation}
To show that $  r^{**}_{_M} $ is the best possible, we consider $ f=f_M $ given by \eqref{e-3.7}. A simple computation, using \eqref{e-3.8}, \eqref{e-3.14} and \eqref{e-3.16} for $ f=f_M $ and $ r>r^{**}_{_M} $ shows that
\begin{align*}
	&\sum_{n=2}^{\infty}\left(|a_n|+|b_n|\right)|z|^n+\sqrt{|J_{f_{M}}(z)}|z|\\&\nonumber= r+\sum_{n=2}^{\infty}\frac{2Mr^n}{n(n-1)}+\left(1+2M\sum_{n=1}^{\infty}\frac{r^n}{n}\right)r\\&>r^{**}_{_M}+\sum_{n=2}^{\infty}\frac{2M\left(r^{**}_{_M}\right)^n}{n(n-1)}+\left(1+2M\sum_{n=1}^{\infty}\frac{\left(r^{**}_{_M}\right)^n}{n}\right)r\\&=2r^{**}_{_M}+2M\left(r^{**}_{_M}+(1-2r^{**}_{_M})\log(1-r)\right)\\&=1+2M(1-2\log2)\\&=d(f_M(0),\partial f_M(\mathbb{D}))
\end{align*}
and hence $ r^{**}_{_M} $ is the best possible radius. This completes the proof.
\end{pf}
\begin{pf}[\bf Proof of Theorem \ref{th-2.7}]
	Let $f \in \mathcal{P}^{0}_{\mathcal{H}}(M)$ be given by \eqref{e-1.2}. In view of Lemma \ref{lem-1.8}, Lemma \ref{lem-1.9}, for $ |z|=r $, we obtain 
	\begin{align}\label{e-33.21}
		& |f(z)|^2+\sum_{n=2}^{\infty}(|a_n|+|b_n|)|z|^n+\frac{r}{1-r}\sum_{n=2}^{\infty}(|a_n|+|b_n|)^2|z|^{2n}\\&\nonumber\leq \left(r+2M\sum_{n=2}^{\infty}\frac{r^n}{n(n-1)}\right)^2+2M\sum_{n=2}^{\infty}\frac{r^n}{n(n-1)}+\frac{4M^2r}{1-r}\sum_{n=2}^{\infty}\frac{r^{2n}}{n^2(n-1)^2}\\&=\nonumber\left(r+2M\sum_{n=2}^{\infty}\frac{r^n}{n(n-1)}\right)^2+2M\sum_{n=2}^{\infty}\frac{r^n}{n(n-1)}\\&\nonumber\quad\quad+\frac{4M^2r}{1-r}\left(\displaystyle\sum_{n=2}^{\infty}\frac{r^{2n}}{n^2}+	\displaystyle\sum_{n=2}^{\infty}\frac{r^{2n}}{(n-1)^2}+\displaystyle\sum_{n=2}^{\infty}-\frac{2r^{2n}}{n-1}+\displaystyle\sum_{n=2}^{\infty}\frac{2r^{2n}}{n}\right).
	\end{align}
	A simple computation shows that
	\begin{equation}\label{e-33.22}
	\begin{cases}
	\displaystyle\sum_{n=2}^{\infty}\frac{r^{2n}}{n^2}=-r^2+{\rm Li}_2\left(r^2\right),\vspace{1.5mm}\\
	\displaystyle\sum_{n=2}^{\infty}\frac{r^{2n}}{(n-1)^2}=r^2{\rm Li}_2\left(r^2\right),\vspace{1.5mm}\\\displaystyle\sum_{n=2}^{\infty}-\frac{2r^{2n}}{n-1}=2r^2\log\left(1-r^2\right),\vspace{1.5mm}\\
	\displaystyle\sum_{n=2}^{\infty}\frac{2r^{2n}}{n}=-2r^2-2\log\left(1-r^2\right).
	\end{cases}
	\end{equation}
Therefore, using \eqref{e-3.5}, \eqref{e-33.21} and \eqref{e-33.22}, we obtain 
\begin{align*}
	&|f(z)|^2+\sum_{n=2}^{\infty}(|a_n|+|b_n|)|z|^n+\frac{r}{1-r}\sum_{n=2}^{\infty}(|a_n|+|b_n|)^2|z|^{2n}\\&\leq \left(r+2M(r+(1-r)\log(1-r))\right)^2+2M\left(r+(1-r)\log(1-r)\right)\\&\quad\quad+\frac{4M^2r}{1-r}\bigg(-3r^2-2\left(1-r^2\right)\log\left(1-r^2\right)+\left(1+r^2\right){\rm Li}_2\left(r^2\right)\bigg).
\end{align*}
It is easy to see that
\begin{align*}
	&\left(r+2M(r+(1-r)\log(1-r))\right)^2+2M\left(r+(1-r)\log(1-r)\right)\\&\quad\quad+\frac{4M^2r}{1-r}\bigg(-3r^2-2\left(1-r^2\right)\log\left(1-r^2\right)+\left(1+r^2\right){\rm Li}_2\left(r^2\right)\bigg)\leq 1+2(1-2\log 2)
\end{align*}
for $ r\leq r^{\prime}_{_M} $, where $ r^{\prime}_{_M} $ is a root of $ H_4(r)=0 $, where $ H_4 : [0,1]\rightarrow\mathbb{R} $ is defined by
\begin{align*}
H_4(r):&=\left(r+2M(r+(1-r)\log(1-r))\right)^2+2M\left(r+(1-r)\log(1-r)\right)\\&\quad\quad+\frac{4M^2r}{1-r}\bigg(-3r^2-2\left(1-r^2\right)\log\left(1-r^2\right)+\left(1+r^2\right){\rm Li}_2\left(r^2\right)\bigg)\\&\quad\quad- 1-2(1-2\log 2)\\&=\left(r+2M\sum_{n=2}^{\infty}\frac{r^n}{n(n-1)}\right)^2+2M\sum_{n=2}^{\infty}\frac{r^n}{n(n-1)}+\frac{4M^2r}{1-r}\sum_{n=2}^{\infty}\frac{r^{2n}}{n^2(n-1)^2}\\&\quad\quad-1-2M\sum_{n=2}^{\infty}\frac{(-1)^{n-1}}{n(n-1)}.
\end{align*}
By the similar argument that being used in the proof of the preceeding results, a simple computation shows that $ r^{\prime}_{_M} $ is the unique root of $ H_4(r) $ in $ (0,1) $, and hence we have
\begin{align}\label{e-33.23}
	&\left(r^{\prime}_{_M}+2M(r^{\prime}_{_M}+(1-r^{\prime}_{_M})\log(1-r^{\prime}_{_M}))\right)^2+2M\left(r^{\prime}_{_M}+(1-r^{\prime}_{_M})\log(1-r^{\prime}_{_M})\right)\\&\nonumber\quad\quad+\frac{4M^2r^{\prime}_{_M}}{1-r^{\prime}_{_M}}\bigg(-3\left(r^{\prime}_{_M}\right)^2-2\left(1-\left(r^{\prime}_{_M}\right)^2\right)\log\left(1-\left(r^{\prime}_{_M}\right)^2\right)\\&\nonumber\quad\quad+\left(1+\left(r^{\prime}_{_M}\right)^2\right){\rm Li}_2\left(\left(r^{\prime}_{_M}\right)^2\right)\bigg)\\&\nonumber= 1+2(1-2\log 2).
\end{align}
To show that $  r^{\prime}_{_M} $ is the best possible radius, we consider $ f=f_M $ given by \eqref{e-3.7}. Then a simple computation, using \eqref{e-3.8}, \eqref{e-3.21} and \eqref{e-33.23} for $ f=f_M $ and $ r>r^{\prime}_{_M} $ shows that
\begin{align*}
& |f_M(z)|^2+\sum_{n=2}^{\infty}(|a_n|+|b_n|)|z|^n+\frac{r}{1-r}\sum_{n=2}^{\infty}(|a_n|+|b_n|)^2|z|^{2n}\\&\nonumber= \left(r+2M\sum_{n=2}^{\infty}\frac{r^n}{n(n-1)}\right)^2+2M\sum_{n=2}^{\infty}\frac{r^n}{n(n-1)}+\frac{4M^2r}{1-r}\sum_{n=2}^{\infty}\frac{r^{2n}}{n^2(n-1)^2}\\&>\left(r^{\prime}_{_M}+2M\sum_{n=2}^{\infty}\frac{\left(r^{\prime}_{_M}\right)^n}{n(n-1)}\right)^2+2M\sum_{n=2}^{\infty}\frac{\left(r^{\prime}_{_M}\right)^n}{n(n-1)}+\frac{4M^2r^{\prime}_{_M}}{1-r^{\prime}_{_M}}\sum_{n=2}^{\infty}\frac{\left(r^{\prime}_{_M}\right)^{2n}}{n^2(n-1)^2}\\&=	\left(r^{\prime}_{_M}+2M(r^{\prime}_{_M}+(1-r^{\prime}_{_M})\log(1-r^{\prime}_{_M}))\right)^2+2M\left(r^{\prime}_{_M}+(1-r^{\prime}_{_M})\log(1-r^{\prime}_{_M})\right)\\&\nonumber\quad\quad+\frac{4M^2r^{\prime}_{_M}}{1-r^{\prime}_{_M}}\bigg(-3\left(r^{\prime}_{_M}\right)^2-2\left(1-\left(r^{\prime}_{_M}\right)^2\right)\log\left(1-\left(r^{\prime}_{_M}\right)^2\right)\\&\nonumber\quad\quad+\left(1+\left(r^{\prime}_{_M}\right)^2\right){\rm Li}_2\left(\left(r^{\prime}_{_M}\right)^2\right)\bigg)\\&\nonumber= 1+2(1-2\log 2)\\&=d(f_M(0),\partial f_M(\mathbb{D}))
\end{align*}
and hence $ r^{\prime}_{_M} $ is the best possible. This completes the proof.
\end{pf}
	
\noindent\textbf{Acknowledgment:} The first author is supported by the Institute Post Doctoral Fellowship of IIT Bhubaneswar, India, the second author is supported by SERB-CRG, India.

\end{document}